\documentclass[11pt]{article}

\usepackage[margin= 2 cm, bottom=18mm,footskip=10mm,top=18mm]{geometry}

\usepackage{amsfonts,amssymb,amsmath,amsthm,amstext}
\usepackage{enumitem}
\usepackage{cleveref}

\newtheorem{theorem}{Theorem}

\theoremstyle{definition}

\usepackage{thmtools,thm-restate}
\usepackage[title]{appendix}

\newcommand{\xb}{\mathbf{x}}

\binoppenalty=\maxdimen
\relpenalty=\maxdimen

\title{A generalization of Tur\'{a}n's theorem}
\author{Domagoj Brada\v{c} \thanks{Department of Mathematics, ETH, Z\"urich, Switzerland. Research supported in part by SNSF grant 200021\_196965. Email: \textbf{domagoj.bradac@math.ethz.ch}.}}
\date{}

\begin{document}
    \maketitle
    \begin{abstract}
        We prove a generalization of Tur\'{a}n's theorem proposed by Balogh and Lidick\'{y}.
    \end{abstract}

    In this short note, we prove the following statement proposed by Balogh and Lidick\'{y} and communicated by Balogh at the Combinatorics, Probability and Computing conference held in Oberwolfach in 2022.
    
    \begin{theorem} \label{thm:thm}
        Let $G$ be an arbitrary graph on $n$ vertices. For each edge $e \in E(G),$ define its weight $w(e)$ as $w(e) = \frac{r}{2(r-1)},$ where $r$ is the size of the largest clique in $G$ containing $e.$ Then, $$\sum_{e \in E(G)} w(e) \le n^2 / 4.$$
    \end{theorem}
    
    In the case of $K_5$-free graphs, \Cref{thm:thm}, along with a corresponding stability result, was proved by Balogh and Lidick\'{y} using flag algebras \cite{balogh_lidicky}.
    
    Observe that \Cref{thm:thm} can be viewed as a generalization of Tur\'{a}n's theorem. Indeed, let $G$ be a $K_{r+1}$-free graph. Since $\frac{t}{2(t-1)} = \frac{1}{2} + \frac{1}{2(t-1)}$ is a decreasing function of $t,$ we have that $w(e) \ge \frac{r}{2(r-1)}$ for every edge $e \in E(G).$ Then, \Cref{thm:thm} implies $e(G) \le \frac{n^2}{4} / \frac{r}{2(r-1)} = \left(1 - \frac{1}{r} \right) \frac{n^2}{2}.$ It is easy to check that \Cref{thm:thm} is tight when $G$ is the Tur\'{a}n graph $T_{n, r}$ and $r$ divides $n.$
    
    Our proof is an adaptation of the proof of Tur\'{a}n's theorem by  Motzkin and Straus \cite{motzkin_straus_1965} using the graph Lagrangian.
    \begin{proof}[Proof of \Cref{thm:thm}]
        Let $G = (V, E)$ with $V = \{v_1, \dots, v_n\}.$ Let $S$ denote the regular $(n-1)$-dimensional simplex, that is, the set of points $(x_1, \dots, x_n) \in \mathbb{R}^n$ such that $x_i \ge 0, \forall i \in [n]$ and $\sum_{i=1}^n x_i = 1.$ For $\xb = (x_1, \dots, x_n) \in S,$ we define 
        \[ f(\xb) = \sum_{v_iv_j \in E} w(v_iv_j) x_i x_j. \]
        Let $m$ be the maximum value of $f$ over $S.$ First we show that $f$ is maximized at some point $\xb = (x_1, \dots, x_n)$ supported on a clique of $G,$ that is, such that $v_iv_j \in E$ whenever $x_i, x_j > 0.$
        
        For the sake of contradiction, suppose there exists $\xb = (x_1, \dots, x_n) \in S$ satisfying $f(\xb) = m$ and $\xb$ has minimum support among all maximums of $f$ (i.e. the minimum number of positive entries), but there exist $i \neq j$ such that $v_iv_j \not\in E$ and $x_i, x_j > 0.$ Define
        \begin{align*}
            &s_i = \sum_{k \colon v_iv_k \in E} w(v_iv_k) \cdot x_k \text{ and}\\
            &s_j = \sum_{k \colon v_jv_k \in E} w(v_jv_k) \cdot x_k.
        \end{align*}
        
        Without loss of generality, assume that $s_i \ge s_j.$ Define $\xb' = (x_1', \dots, x_n')$ as follows:
        \[ x_k' = 
        \begin{cases} 
                x_k, &\text{ if } k \not\in \{i, j\},\\
                x_i + x_j, &\text{ if } k = i,\\
                0, &\text{ if } k = j.
        \end{cases} \]
        Clearly $\xb' \in S$ and observe that
        \begin{align*}
            f(\xb') - f(\xb) &= \sum_{k \colon v_iv_k \in E} w(v_iv_k) \cdot (x_i' - x_i) x_k + \sum_{k \colon v_jv_k \in E} w(v_jv_k) \cdot (x_j' - x_j)x_k \\ &= \sum_{k \colon v_iv_k \in E} w(v_iv_k) \cdot x_j x_k + \sum_{k \colon v_jv_k \in E} w(v_jv_k)\cdot  (-x_j) x_k = x_j (s_i - s_j) \ge 0.
        \end{align*}
        Hence, $\xb'$ is also a maximum of $f,$ but its support is stricly smaller, contradicting our assumption on $\xb.$
        
        Now, let $\xb$ be a maximum of $f$ supported on a clique $K$ of size $r$ and without loss of generality assume that $K = \{ v_1, \dots, v_r\}.$ As already noted, $\frac{t}{2(t-1)}$ is a decreasing function of $t,$ so $w(v_iv_j) \le \frac{r}{2(r-1)}$ for all $1 \le i < j \le r.$ Therefore,
        \begin{align*}
            f(\xb) &= \sum_{1 \le i < j \le r} w(v_iv_j) \cdot x_ix_j \le \frac{r}{2(r-1)} \sum_{1 \le i < j \le r} x_ix_j = \frac{r}{2(r-1)} \cdot \frac{1}{2} \sum_{1 \le i \le r} \sum_{\substack{1 \le j \le r\\ j\neq i}} x_ix_j \\
            &= \frac{r}{4(r-1)} \sum_{1 \le i \le r} x_i(1 - x_i) = \frac{r}{4(r-1)} \left(1 - \sum_{1 \le i \le r} x_i^2 \right) \le \frac{r}{4(r-1)} \left( 1 - r \cdot \frac{1}{r^2} \right) = \frac{1}{4},
        \end{align*}
        where in the third and fourth equalities we used that $\sum_{i=1}^r x_i = 1$ and in the last inequality we used Jensen's inequality. We conclude that $m \le \frac{1}{4}.$ On the other hand, setting $\xb = \left( \frac{1}{n}, \dots, \frac{1}{n}\right),$ we obtain
        \[ \frac{1}{4} \ge m \ge f(\xb) = \sum_{v_iv_j \in E} w(v_iv_j) / n^2, \]
        which, after rearranging, finishes the proof.
    \end{proof}
    
    \textbf{Acknowledgment.}
    The author would like to thank Benny Sudakov for introducing him to the problem and useful discussions.
    
    \bibliographystyle{abbrv}
    
\end{document}